\documentclass[reqno,12pt,draft]{amsart}

\usepackage{amssymb,amsmath,amsthm}

     \makeatletter
\def\multilimits@{\bgroup
  \Let@
  \restore@math@cr
  \default@tag
 \baselineskip\fontdimen10 \scriptfont\tw@
 \advance\baselineskip\fontdimen12 \scriptfont\tw@
 \lineskip\thr@@\fontdimen8 \scriptfont\thr@@
 \lineskiplimit\lineskip
 \vbox\bgroup\ialign\bgroup\hfil$\m@th\scriptstyle{##}$\hfil\crcr}
\def\Sb{_\multilimits@}
\def\Sp{^\multilimits@}
\def\endSb{\crcr\egroup\egroup\egroup}

     \makeatother

\normalbaselines

\newtheorem{thm}{\hspace{\parindent}Theorem}
\newtheorem{lem}[thm]{\hspace{\parindent}Lemma}

\newtheorem{defn}[thm]{\hspace{\parindent}Definition}
{\theoremstyle{remark}

}

\newcommand{\ip}[1]{\langle #1 \rangle}
\newcommand{\nrm}[1]{|\!| #1 |\!|}
\newcommand{\ci}[1]{_{{}_{\scriptstyle{#1}}}}
\newcommand{\R}{\mathbb{R}}
\newcommand{\T}{\mathbb{T}}

\newcommand{\f}{\varphi}
\newcommand{\al}{\alpha}

\begin{document}

\title[Eventual Regularity of the Solutions to the Supercritical SQG]{Eventual Regularity of the Solutions to the Supercritical Dissipative Quasi-Geostrophic Equation}
\author{Michael Dabkowski}
\begin{abstract} Recently in  \cite{S}, Silvestre proved that certain weak solutions of the slightly supercritical surface quasi-geostrophic equation eventually become smooth. To prove this, he employed a De Giorgi type argument originated in the work of Caffarelli and Vasseur, \cite{CV}. In \cite{KN}, Kiselev and Nazarov  proved a variation of the result of Caffarelli and Vasseur by introducing a class of test functions. Motivated by the results of Silvestre, we will modify the class of test functions from \cite{KN} and use this modified class to show that a solution to the supercritical SQG that is smooth up to a certain time must remain smooth forever.     
\end{abstract}

\maketitle

\section{Introduction}

The setting of this paper will be the $d$-dimensional torus, $\T^d$. We may equivalently think of the problem in the setting of $\R^d$ with periodic initial data. Throughout the paper we will consider only real valued functions. We consider the Cauchy problem for the dissipative equation 
\begin{equation}
\left\{ \begin{array}{l}
\theta_t = (u \cdot \nabla) \theta - (-\Delta)^{\alpha/2} \theta\\
\theta(x,0)= \theta_0(x)\\
\end{array}\right., \label{SQG}
\end{equation} 
where $u = R\theta$, $R$ is a certain divergence free operator, and $(-\Delta)^{\alpha/2}$ is the fractional Laplacian. In the case of the surface quasi-geostrophic equation (SQG for brevity), $d=2$ and $u = (-R_2\theta, R_1\theta)$, where the $R_j$s are the standard Riesz transforms. These operators are defined on a suitably smooth class of functions by multiplication on the Fourier side. For $n \in \mathbb{Z}^d$, if 
$$\widehat{\theta}(n) = \frac{1}{(2\pi)^d}\int_{\T^d} \theta(x) e^{-i n \cdot x}\,dx $$
is the $n^{th}$ Fourier coefficient of $\theta$, then for $n\neq 0$
$$ \widehat{(R_j\theta)}(n) = i\frac{n_j}{|n|} \widehat{\theta}(n) \quad \mbox{and} \quad [(-\Delta)^{\alpha/2}\theta]\,\,\widehat{}\,\,(n) = |n|^{\alpha} \widehat{\theta}(n),$$
and $\widehat{(R_j\theta)}(0) = [(-\Delta)^{\alpha/2}\theta]\,\,\widehat{}\,\,(0)= 0.$ 

The parameter $\al$ ranges between $0$ and $2$. The case when $\al \in (1,2]$ is referred to as the subcritical case. In the subcritical case, the global well-posedness has been established in the case of smooth initial data (See \cite{CCW} and the references therein). The critical case, $\al=1$, has been the source of much study in recent years. In \cite{CCW}, Constantin, Cordoba, and Wu proved that if the $L^{\infty}$ norm on the initial data is small enough, then there is a global regular solution. Later, Kiselev, Nazarov, and Volberg introduced the modulus of continuity method in \cite{KNV}. This method was used to prove the global well posedness of the critical SQG for smooth periodic initial data by finding a priori bounds on $\nrm{\nabla \theta}_{\infty}$. In the supercritical case, $\al <1$, many open questions remain. 

In \cite{CV}, Caffarelli and Vasseur used De Giorgi iteration to show that a uniform bound in $BMO$ of the velocity term in the drift diffusion equation implies that certain weak solutions are locally H\"{o}lder continuous.  In \cite{KN}, Kiselev and Nazarov showed that, in the case of the critical surface quasi-geostrophic equation, a uniform bound in $BMO$ on a smooth velocity leads to a certain degree of H\"{o}lder continuity. In this way they were able to give yet another proof of the existence of global smooth solutions to the critical surface quasi-geostrophic equation. Their approach relies on passing the evolution onto a special class of functions which is ``dual" to the class of  H\"{o}lder continuous functions. The reason for studying the H\"{o}lder continuity of solutions to the SQG can be seen from results of Constantin and Wu \cite{CW}. They showed that if you have a uniform bound on the $C^{1- \al + \delta}$ ($\delta >0$) norm of a certain weak solution to the SQG on a time interval, then in fact you have a smooth classical solution on that interval. Recently in \cite{S}, Silvestre proved that if the dissipative power is slightly smaller than $1/2$, namely the power of the Laplacian is $\frac{1-\epsilon}{2}$ for small $\epsilon$, then certain weak solutions become H\"{o}lder continuous after a certain time. The proof employs De Giorgi-type estimates to show that on a parabolic cylinder the oscillation of a certain continuation of the solution is not more than a fraction of the oscillation of the continuation of the solution on a twice larger parabolic cylinder under the assumption that the $L^{2d/\al}$ norm of the velocity is uniformly bounded. 

Before stating the main result, we define  
$$ \nrm{f}_{k,p} = \nrm{\nabla^k f}_{p},$$
which for mean zero functions can be shown to be equivalent to the standard Sobolev norm by the Poincar\'{e} inequality. Motivated by the smooth class constructed in the work of Nazarov and Kiselev and the work of Silvestre, we will prove the following theorem:

\begin{thm} Suppose that $R$ is a divergence free vector-valued operator that, for every $k\geq 0$ and every $1<p<\infty$, satisifes $\nrm{Rf - Rg}_{k,p} \leq C(k,p) \nrm{f-g}_{k,p}$ for some constants $C(k,p)$, and, for every $\epsilon >0$ satisifes $\nrm{\nabla (R f)}_{\infty} \leq C(\epsilon) \nrm{\nabla f}_{C^{\epsilon}}$ for some constant $C(\epsilon)$. There is a time $T = T(\al,\nrm{\theta_0}_{\infty})$ such that if $\theta \in C^{\infty}(\T^d \times [0,T])$ is a solution to the Cauchy problem
$$\left\{ \begin{array}{l}
\theta_t = (R \theta \cdot \nabla) \theta - (-\Delta)^{\alpha/2} \theta\\
\theta(x,0)= \theta_0(x)\\
\end{array}\right.,
$$
then $\theta$ extends to a solution in $C^{\infty}(\T^d \times [0,\infty)).$
\end{thm}

A consequence of this theorem is 

\begin{thm}[Eventual Regularization for the Supercritical SQG] There is a time $T = T(\al,\nrm{\theta_0}_{\infty})$ such that if $\theta \in C^{\infty}(\T^2 \times [0,T])$ is a solution to
$$
\left\{ \begin{array}{l}
\theta_t = (R^{\perp}\theta \cdot \nabla) \theta - (-\Delta)^{\alpha/2} \theta\\
\theta(x,0)= \theta_0(x)\\
\end{array}\right.,
$$
then $\theta$ extends to a solution in  $C^{\infty}(\T^2 \times [0,\infty)).$ 
\end{thm}

Classical results about Riesz transforms imply $R^{\perp}$ satisfies the conditions in Theorem $1$ (See \cite{St}). Both Theorems tell us that for any value of $\al$ in the supercritical range, if we have a solution that is smooth up to a certain time, then it remains smooth forever.

\section{Dualizing the Problem}

We now define a variant of the class introduced in \cite{KN}. Let $A > 1$ be a parameter to be fixed later. 
\begin{defn} We will say that a smooth function $\psi$ defined on $\T^d$ is in $\mathcal{U}(r)$ if 
\begin{equation}
\int_{\T^d} |\psi(x)|^p \,dx \leq A r^{-(p-1)d} \label{size}
\end{equation}
and
\begin{equation}
\sup\Big\{ \Big| \int_{\T^d} f(x)\psi(x)\,dx \Big| : f \in C^{\infty} \cap \text{Lip}(1)\Big\} \leq r \label{smooth}
\end{equation}

\end{defn}
In this definition we have used the notation $Lip(M)$ to denote the class of all functions $f$ such that $|f(x) -f(y)| \leq M|x-y|$ for all $x,y\in \T^d.$
Since all constant functions are Lipschitz, if $(\ref{smooth})$ holds, then the function $\psi$ must have mean zero. Also notice that if a function $\psi$ satisfies $|\!|\psi|\!|_1 \leq 1$ and $|\!|\psi|\!|_{\infty} \leq A^{1/(p-1)} r^{-d}$ (as were the conditions in \cite{KN}), then interpolation shows that $\psi$ satisfies $(\ref{size})$. If $\f$ is supported in $B_r$ (the ball of radius $r$ centered at the origin in $\T^d$), has mean zero and $\nrm{\f}_{p} \leq r^{-d/q}$, then $\f \in \mathcal{U}(r)$. In what follows we will write $f \in B\,\mathcal{U}(r)$ if $f/B \in \mathcal{U}(r).$ 

Recall that a function $g$ is H\"{o}lder continuous with exponent $\beta \in (0,1)$ if $|g(x)-g(y)| \leq C |x-y|^{\beta}$ for some constant $C$ and all $x,y \in \T^d.$ We will denote the class of  H\"{o}lder continuous functions with exponent $\beta$ on $\T^d$ by $C^{\beta}(\T^d)$ . Paley-Littlewood projections can be used to characterize $C^{\beta}(\T^d)$ as follows: we let $\omega$ be a smooth compactly supported function on $\mathbb{R}^d$ that is identically $1$ when $|x| \leq 1$, radially decreasing, and vanishing for $|x|\geq 2$. Define $\f(x) = \omega(x)- \omega(2x)$ and $\f_j(x) = \f(x/2^j).$ For an integrable function $f$ on $\T^d$, we define for any non-negative integer $j$
$$ \Delta_j f (x) =  \sum_{n \in \mathbb{Z}^d} \f_j (n) \widehat{f}(n) e^{in\cdot x}.$$
The operators $\Delta_j$ are essentially smooth projections on the frequency scale $2^{j}$.  Recall that a bounded function $g$ on $\T^d$ is $C^{\beta}(\T^d)$ if and only if for every $j \geq 0$,
\begin{equation}
\nrm{\Delta_j g}_{\infty} \leq M 2^{-\beta j} ,\label{CB}
\end{equation}
with some $M>0$. See \cite{K} for a proof. 
Write $\ip{f,g}$ to denote $\int_{\T^d} f(x)g(x)\,dx.$ If we have control over 
$$|\ip{g,\mathcal{U}(r)}| := \sup \big\{|\ip{g,\psi}|: \psi \in \mathcal{U}(r) \big\},$$ 
then we get control over the H\"{o}lder $C^{\beta}$ seminorm of $g$. More precisely, we have the following 
\begin{lem} Suppose that a function $g$ on $\T^d$ has the property that 
$$ \Big| \int_{\T^d} g(x) \psi(x) \,dx \Big| \leq r^{\beta},$$
for all $\psi \in \mathcal{U}(r)$ and $0<r\leq 1$. Then $g \in C^{\beta}(\T^d)$ and 
$$ \nrm{g}_{C^{\beta}} := \sup_{x\neq y} \frac{|g(x) - g(y)|}{|x-y|^{\beta}} \leq C(\beta).$$
\end{lem}

\proof Let $\f_j$ be the function defined above and let $\mathcal{F}^{-1}$ denote the inverse Fourier transform. Since $\f$ is smooth, compactly supported, and vanishes at the origin, $\mathcal{F}^{-1}{\f}$ is a Schwartz function with mean zero. We can therefore find a constant $C$ such that $\int_{\R^d} |\mathcal{F}^{-1}{\f}(x)| \,dx \leq C, \int_{\R^d}|x|\cdot|\mathcal{F}^{-1}{\f}(x)| \,dx \leq C$, and $|\mathcal{F}^{-1}\f(y)| \leq C (1+|y|^2)^{-d}$ for all $y \in \R^d$. Scaling these inequalities gives for all $j\geq 0$,
\begin{enumerate}
\item[(i)] $\int_{\R^d} |\mathcal{F}^{-1}{\f_j}(x)| \,dx \leq C$
\item[(ii)] $\int_{\R^d}|x|\cdot|\mathcal{F}^{-1}{\f_j}(x)| \,dx \leq C2^{-j}$
\item[(iii)] $|\mathcal{F}^{-1}\f_j(y)| \leq C 2^{jd}(1+|2^jy|^2)^{-d}$ for all $y \in \R^d$
\end{enumerate}

Now define for $x \in \T^d,$  
$$\Phi_j(x) = c \sum_{n \in \mathbb{Z}^d} (\mathcal{F}^{-1} \f_j) (x+2\pi n),$$
for some constant $c$ to be chosen later.  We claim that if we choose $c$ sufficiently small independently of $j$, then $\Phi_j \in \mathcal{U}(2^{-j}).$ Inequality (iii) implies that $\nrm{\Phi_j}_{\infty} \leq cC' 2^{jd}$ and inequality (i) implies $\nrm{\Phi_j}_{1} \leq cC'$ for some constant $C'>C$. Interpolation shows the norm condition, $(\ref{size})$, is satisfied provided $c$ is small. Inequality (ii) and mean zero property of $\mathcal{F}^{-1}\f_j$ show that
$$ \Big| \int_{\R^d} f(x) \mathcal{F}^{-1} \f_j (x) \,dx \Big| \leq C 2^{-j},$$
for any smooth $2\pi$-periodic function $f$ in $Lip(1)$ on $\R^d$. This implies that the smooth condition, $(\ref{smooth})$, is satisfied, provided $c$ is small enough. 
Since $\Delta_j$ is a convolution with $\Phi_j$ on the space side, we have for any $y \in \T^d$,
$$ \int_{\T^d} g(x) \Phi_j(x-y)\,dx = c \Delta_j g(y).$$
Since the class $\mathcal{U}(r)$ is invariant under translations $\Phi_j(\cdot-y) \in \mathcal{U}(2^{-j})$ for any $y\in \T^d$, which implies the left hand side of the above equality is not more than $2^{-j\beta}$ by assumption. Therefore, we have $\nrm{\Delta_j g}_{\infty} \leq c^{-1} 2^{-j\beta}$ for all $j\geq 0$. From $(\ref{CB})$ we conclude that $g \in C^{\beta}(\T^d)$. \qed 
\vspace{.1in}

If we wish to prove that the solution to the SQG at time $t$ is in $C^{\beta}$, we must estimate $|\ip{\theta(\cdot, t), \mathcal{U}(r)}|$. We do so by determining how the class $\mathcal{U}(r)$ evolves under the backward equation. Let $u = R\theta$ and suppose that $\psi$ is a solution to the equation
\begin{equation}
\psi_s = (u \cdot \nabla) \psi + (-\Delta)^{\alpha/2} \psi \label{URSQG}
\end{equation}
Later we will impose some ``future condition" that at some moment of time, the solution to $(\ref{URSQG})$ is in the class $\mathcal{U}(r)$. Consider the pairing function defined for $\tau>0$,
$$P(\tau)=\ip{\theta(\cdot, \tau), \psi(\cdot, \tau)} = \int_{\T^d} \theta(x, \tau) \psi(x,\tau)\,dx$$ 
Equations $(\ref{SQG})$ and $(\ref{URSQG})$ together with the fact that $u$ is a divergence free vector field imply that $P$ has zero derivative and is therefore constant:
\begin{equation}
\int_{\T^d} \theta(x, t-s) \psi(x,t-s)\,dx = \int_{\T^d} \theta(x,t) \psi(x,t)\,dx, \label{duality}
\end{equation}
for any pair of times $t\geq s\geq 0.$ 
Our next objective will be to find functions $F$ and $G$ such that if $\psi$ is a solution to $(\ref{URSQG})$ and $\psi(\cdot,t) \in \mathcal{U}(r)$ for some fixed time $t$, then $\psi(\cdot, t-s) \in F(s,r)\mathcal{U}(G(s,r))$ for $s$ sufficiently small. This is what is meant by dualizing the problem: $(\ref{duality})$ allows us to move the dynamics from a solution to $(\ref{SQG})$ onto a solution of $(\ref{URSQG})$. Now we can determine the H\"{o}lder regularity of the solution to $(\ref{SQG})$ by determining how $(\ref{URSQG})$ alters the class $\mathcal{U}(r)$.  

\section{Evolution of the Class $\mathcal{U}(r)$}

The heart of the proof of Theorem $1$ is the following

\begin{lem} (Class Evolution) Given $\al,\beta \in (0,1)$ and $p>1$ such that $\al +\beta -d/q >1$ where $q$ is the conjugate exponent to $p$, there are parameters $\delta, r_0>0$ with the following property: if $0< r \leq r_0$, $s \leq r^{\al}$, and $|\ip{\theta(\cdot, \tau),\mathcal{U}(R)}| \leq R^{\beta}$ for all $R\geq r e^{\delta} $ and $\tau \in [t-s,t]$, then every solution $\psi$ to $(\ref{URSQG})$ with $\psi(\cdot,t) \in \mathcal{U}(r)$ satisfies
\begin{equation}
\psi(\cdot, t-s) \in \exp(-\delta s r^{-\al})\mathcal{U}\big(r\exp(\delta \beta^{-1} sr^{-\al})\big).
\label{evolve}
\end{equation}
\end{lem}

\proof Let $\widetilde{\chi}$ be a smooth radially decreasing function on $\R^d$ supported in $\{x: |x| \leq 1\}$ and has mean $1$. Then define $\widetilde{\chi}_r(x) = r^{-d} \widetilde{\chi}(x/r)$, so that $ \widetilde{\chi}_r$ has mean $1$ for every $r>0$ and $\nrm{\nabla \widetilde{\chi}_r}_{\infty} \leq Cr^{-d-1}$. For $r\in(0,1]$, the function $\widetilde{\chi}_r$ is supported in $(-\pi,\pi)^d$. Using the identification of $\T^d$ with $[-\pi,\pi]^d$, we identify $\widetilde{\chi}_r$ with a function $\chi_r$ on $\T^d$. 

Suppose now that $f = f(x,\tau)$ is a solution to the smooth forward evolution:
\begin{equation}
\left\{ \begin{array}{l}
f_{\tau}= (u_{r} \cdot \nabla) f - (-\Delta)^{\alpha/2} f\\
f(x,t-s) = f_0(x)
\end{array}\right., \label{LIPSQG}
\end{equation}
where $u_{r} = R(\theta \ast \chi_r)$ and $R$ is the divergence free operator in the statement of Theorem $1$.
Consider the pairing function $\Pi$ defined in the interval $[t-s,t]$ by
$$\Pi(\tau) = \ip{f(\cdot,\tau), \psi(\cdot,\tau)} = \int_{\T^d} f(x,\tau) \psi(x,\tau)\,dx,$$
where $\psi$ is a solution to $(\ref{URSQG})$. Differentiating $\Pi$ and using that $u_{r}$ and $u$ are divergence free, we see that 
$$\Pi'(\tau) = \int_{\T^d} [(u_{r}(x,\tau) - u(x,\tau))\cdot \nabla f(x,\tau)] \psi(x,\tau)\,dx.$$
Integrating the above expression we have
\begin{equation}
\int_{\T^d} \psi(x,t-s) f_0(x) \,dx = \int_{\T^d} \psi(x,t) f(x,t)\,dx + \int_{t-s}^{t} \int_{\T^d} [(u-u_{r})(x,\tau)\cdot \nabla f(x,\tau)] \psi(x,\tau) \,dx\,d\tau. \label{lipduality}
\end{equation}
Call the absolute value of the first integral on the right hand side of the above equation the \emph{smooth part}, which we denote by $I$, and the absolute value of the second integral the \emph{rough part}, which we denote by $II$. Then we have  
$$ \Big|\int_{\T^d} \psi(x,t-s) f_0(x) \,dx \Big| \leq I + II.$$
In what follows we will estimate the quantities $I$ and $II$.

\section{Modulus of Continuity Redux: The Smooth Part}

Using a modulus of continuity argument inspired by \cite{KNV}, we will determine the size of the smooth part.
\begin{lem}(Lipschitz Evolution)
Let $v$ be a smooth divergence free vector field such that $v(x,\tau) \in Lip(M)$ for all $\tau \in [t-s,t].$ Suppose that $f$ is a solution to the system 
$$\left\{ \begin{array}{l}
f_{\tau} = (v \cdot \nabla) f - (-\Delta)^{\alpha/2} f\\
f(x,t-s) = f_0(x)
\end{array}\right.,
$$ 
where $f_0$ is a smooth function in $Lip(1)$. Then $f(x,t-k) \in Lip(\exp(M(s-k))),$ for all $k\in[0,s].$ 
\end{lem}

\proof Fix $\epsilon > 0$ and consider 
$$ \kappa = \sup \big\{k\in[0,s] : \exists x,y \in \T^d \,\,\text{such that}\,\, |f(x,t-k)-f(y,t-k)| \geq e^{M(s-k)}(|x-y| + \epsilon)\big\}.$$
The global regularity theory for an equation of the form $(\ref{LIPSQG})$ with smooth velocity implies that $f$ is smooth for all times. In what follows we will omit the absolute value signs around the quanitity $f(x,\tau)-f(y,\tau)$ as we may always make it is non-negative by exchanging $x$ and $y$ if necessary. Suppose $\kappa \geq 0$. 
First, we notice that $\kappa \neq s$. Indeed, if it were, then there would be sequences of points $x_n,y_n \in \T^d$ and $k_n \rightarrow s$ such that 
$$ f(x_n,t-k_n) - f(y_n, t-k_n) \geq \exp(M(s-k_n))(|x_n-y_n|+\epsilon).$$
The compactness of $\T^{2d}$ implies there are points $x$ and $y$ such that
$$ f(x,t-s) - f(y,t-s) \geq  |x-y|+\epsilon,$$
which contradicts the assumption on $f_0$. 
It follows that for $\kappa < k < s$ and all $x,y \in \T^d$ we have 
$$f(x,t-k)-f(y,t-k) < \exp(M(s-k))(|x-y| + \epsilon).$$
Passing to the limit as $k$ tends to $\kappa$ in the previous inequality, the continuity of $f$ implies that for all $x,y \in \T^d$,
$$f(x,t-\kappa) - f(y,t-\kappa) \leq \exp(M(s-\kappa))(|x-y| + \epsilon).$$ 
Using the same compactness argument as above, we see that there are points $x,y \in \T^d$ such that
\begin{equation}
f(x,t-\kappa) - f(y,t-\kappa) = \exp(M(s-\kappa))(|x-y| + \epsilon).\label{lipequal}
\end{equation}
We now claim that there is a $k > \kappa$ such that at these points $x,y \in \T^d$ we have 
\begin{equation}
f(x,t-k) - f(y,t-k) \geq \exp(M(s-k)(|x-y| + \epsilon).\label{liplarge}
\end{equation}
To this end, we will now compute 
\begin{equation}
\partial_{k} (f(x, t-k) - f(y, t-k))\Big|_{k= \kappa}. \label{derivative}
\end{equation}

The velocity term is the derivative of $f$ in the direction of $v$; more precisely, the chain rule gives $(v \cdot \nabla f)(x) = \frac{d}{dh} f(x + h v(x)) |_{h=0}.$ At the breaking points $x$ and $y$ we have
$$ f(x+ hv(x),t-\kappa) - f(y + hv(y),t-\kappa) \leq \exp(M(s-\kappa))[(|x-y|  + h |v(x) - v(y)|) + \epsilon].$$  
Subtracting $\exp(M(s-\kappa))(|x-y|+\epsilon)$ from both sides, dividing by $h$, and passing to the limit gives
$$  (v \cdot \nabla f)(x,t-\kappa) - (v \cdot \nabla f)(y,t-\kappa) \leq \exp(M(s-\kappa)) M |x-y|.$$ 

The next contribution to $(\ref{derivative})$ comes from the dissipative term. Consider the pure dissipative equation 
$$\left\{ \begin{array}{l}
g_{\tau} = - (-\Delta)^{\alpha/2} g\\
g(\cdot,t-\kappa) = f(\cdot,t-\kappa)
\end{array}\right..
$$ 
The solution to this equation is $g(z,\tau) = f(\cdot,t-\kappa) \ast \Phi(z,\tau)$, where $\widehat{\Phi}(\xi,\tau) = \exp(-|\xi|^{\al} \tau)$. The estimates on $f$ at time $t-\kappa$ imply 
$$g(z_1,t-\kappa) - g(z_2,t-\kappa) \leq \exp(M(s-\kappa))(|z_1-z_2| + \epsilon),$$
for all $z_1,z_2 \in \T^d$. Since the solutions to the purely dissipative equation perserve the modulus of continuity, $g(x,\tau) - g(y,\tau) \leq \exp(M(s-\kappa))(|x-y|+\epsilon)$ for all $\tau \geq t-\kappa$.  The contribution of the dissipative part to $(\ref{derivative})$ is exactly the same as $\partial_{\tau}(g(x,\tau) -g(y,\tau))|_{\tau =t-\kappa}$. Since $g(x,t-\kappa) - g(y,t-\kappa) = \exp(M(s-\kappa))(|x-y|+\epsilon)$ and $g(x,\tau) - g(y,\tau) \leq \exp(M(s-\kappa))(|x-y|+\epsilon)$ for all $\tau \geq t-\kappa$, the function $g(x,\tau) - g(y,\tau)$ has a local maximum at $\tau = t-\kappa$. It follows that $\partial_{\tau}(g(x,\tau) -g(y,\tau))|_{\tau =t-\kappa} \leq 0$ and 
$$ [- (-\Delta)^{\al/2} f] (x,t-\kappa) - [- (-\Delta)^{\al/2} f] (y,t-\kappa) \leq 0.$$ 

Combining the estimates for the velocity and the dissipation we see that
\begin{equation} 
\partial_{k} (f(x,t-k) - f(y,t-k)) \Big|_{k=\kappa} \geq - \exp(M(s-\kappa))M|x-y| \label{lipest}
\end{equation}
The $k$ derivative of the growth condition, $\exp(M(s-k))(|x-y|+ \epsilon)$, at the point $\kappa$ is $-M\exp(M(s-\kappa))(|x-y|+ \epsilon)$, which is strictly smaller than the right hand side of $(\ref{lipest})$. It follows that for $k$ slightly larger than $\kappa$ we have $(\ref{liplarge})$, which contradicts the choice of $\kappa$. From this we conclude
$$\big\{k\in[0,s] : \exists x,y \in \T^d \,\,\text{such that}\,\, |f(x,t-k)-f(y,t-k)| \geq e^{M(s-k)}(|x-y|+\epsilon) \big\}$$
is empty for every $\epsilon>0$ and the lemma follows. \qed 
\vspace{.1in}

We now wish to apply the previous lemma to the solution $f$ of $(\ref{LIPSQG})$. Since $\theta_r=\theta \ast \chi_r$ is a smooth function, the velocity term $u_{r} = R (\theta_r)$ is a $C^1$ divergence free vector field. In order to find a uniform bound on the Lipschitz constant of $u_{r}$, we must first estimate the H\"{o}lder norm of $\nabla\theta_r$. To this end, we notice that the choice of $\chi_r$ implies that  $\nabla \chi_r$ is a mean zero vector-valued function supported in the set $B_r$ and has $L^{\infty}$ norm at most $Cr^{-d-1}$. It follows that $r \nabla \chi_r \in C\mathcal{U}(2r)$ (meaning each component is in $C\mathcal{U}(2r)$). By the assumption of the Class Evolution Lemma, after the time $t-s$ the solution pairs well against $\mathcal{U}(R)$ for $R \geq re^{\delta}$; therefore provided $\delta<1/2$, we see that for any $\tau \in[t-s,t]$ we have  $\nrm{\nabla \theta_r(\cdot,\tau)}_{\infty} = \nrm{\theta(\cdot,\tau) \ast \nabla \chi_r}_{\infty} \leq C|\ip{\theta(\cdot, \tau), r^{-1}\mathcal{U}(2r)}| \leq C r^{\beta -1}.$ Similarly, $\nrm{\nabla(\theta(\cdot,\tau) \ast \nabla \chi_r)}_{\infty} \leq C r^{\beta -2}.$  Let $\epsilon >0$. Interpolation implies that the $C^{\epsilon}$ norm of $\nabla \theta_r$ is no more than $Cr^{\beta -1 - \epsilon}$. The norm assumption on $R$ implies $\nrm{\nabla R \theta_r(\cdot, \tau)}_{\infty} \leq C r^{\beta-1-\epsilon}$. Thus, $R \theta_r(\cdot, \tau) \in Lip(Cr^{\beta -1 -\epsilon})$ for any $\tau \in [t-s,t]$. The Lipschitz Evolution Lemma implies $f(\cdot,t) \in Lip (\exp(Csr^{\beta-1-\epsilon}).$ Since $\psi(\cdot, t) \in \mathcal{U}(r)$ by assumption, we have the following estimate for the smooth part:
\begin{equation}
I \leq \Big| \int_{\T^d} \psi(x,t) f(x,t)\,dx \Big| \leq r\exp(Csr^{\beta-1-\epsilon}).\label{S}
\end{equation} 

\section{Mean Zero Duality: The Rough Part}

In this section, we will estimate the rough part of the evolution. Recall that the rough part was the expression
$$ II = \Big|  \int_{t-s}^{t} \int_{\T^d} \big((u-u_{r})(x,\tau)\cdot \nabla f(x,\tau)\big) \psi(x,\tau) \,dx\,d\tau \Big|.$$ 
Trivially estimating $II$ by H\"{o}lder's inequality yields
\begin{equation}
II \leq s \sup_{\tau \in [t-s,t]} \big(\nrm{(u-u_{r})(\cdot,\tau)}_{q}  \nrm{\psi(\cdot,\tau)}_{p} \nrm{\nabla f(\cdot,\tau)}_{\infty}\big).\label{rough1}
\end{equation}
The maximum principle implies $\nrm{\psi(\cdot,\tau)}_{p} \leq \nrm{\psi(\cdot,t)}_{p} \leq A^{1/p} r^{-d/q}$ for $\tau \in [t-s,t]$. The Lipschitz Evolution Lemma implies that $\nrm{\nabla f(\cdot,\tau)}_{\infty}$ is not more than $\exp(Csr^{\beta-1-\epsilon})$. Since $R$ is Lipchitz in the $L^q$ norm,  $\nrm{(u-u_{r})(\cdot,\tau)}_{q} \leq C_q \nrm{(\theta -\theta_r)(\cdot,\tau)}_{q},$  so it suffices to bound $\nrm{(\theta -\theta_r)(\cdot,\tau)}_{q}.$ 

Recall that $\chi_r$ was chosen to have mean $1$, so for any constant $c$ 
\begin{equation}
\nrm{(\theta-\theta_r)(\cdot, \tau) \chi \ci{B_{r}}}_q \leq \nrm{(\theta-c)(\cdot, \tau) \chi \ci{B_{r}}}_q + \nrm{(c-\theta)_r(\cdot, \tau) \chi \ci{B_{r}}}_q \leq 2 \nrm{(\theta-c)(\cdot, \tau) \chi \ci{B_{3r}}}_q. \label{uncertain}
\end{equation}
We claim that for some choice of $c$ the above expression is not more than $Cr^{\beta +d/q}.$ We will prove this using the smoothness on larger scales along with the following 
 
\begin{lem}[Mean Zero Duality] For any $\rho >0,$ there is a constant $c$ such that
\begin{equation} 
\Big(\int_{B_{\rho}} |\theta-c|^{q}\,dx\Big)^{1/q} \leq \sup \{\rho^{d/q} |\ip{\theta,\psi}| : \psi \in \mathcal{U}(\rho)\}.\label{meanzeroduality}
\end{equation} 
\end{lem}

\proof Choose $c$ so that  $\text{sgn}(\theta-c)|\theta-c|^{q-1}$ has mean zero on $B_{\rho}$ and define $\lambda = \rho^{-d/q} \nrm{(\theta-c)\chi \ci{B_{\rho}}}_{q}^{-q/p}.$ With these choices $\psi = \lambda\, \text{sgn}(\theta-c)|\theta-c|^{q-1} \chi \ci{B_{\rho}}$ is a mean zero function supported in $B_{\rho}$. A direct computation shows $\nrm{\psi}_p^p \leq  \rho^{-(p-1)d} \leq A \rho^{-(p-1)d}$ since $A >1$. As mentioned previously, the mean zero, support, and norm properties of $\psi$ imply the Lipshitz pairing condition. We know choose a sequence of smooth functions $\psi_j \in \mathcal{U}(\rho)$ which converge to $\psi$ in $L^p$ norm. Since $\rho^{d/q} |\ip{\theta,\psi_j}| \rightarrow \rho^{d/q} |\ip{\theta,\psi}|$ and the latter expression is left hand side of $(\ref{meanzeroduality})$, the Lemma follows. \qed 
 
\vspace{.1in}

Applying the lemma to the left hand side of $(\ref{uncertain})$ with $\rho = 3r$ for any $\tau \in [t-s,t]$ gives
$$ \Big(\int_{B_r} |(\theta- \theta_r)(x,\tau)|^q \,dx\Big)^{1/q} \leq C(3r)^{d/q} |\ip{\theta, \mathcal{U}(3r)}| \leq C r^{d/q+\beta}.$$
Since $\T^d$ can be covered by a constant multiple of $r^{-d}$ balls of radius $r$, adding the $q^{th}$ powers of the left hand sides of the above inequalities for all these balls yields 
\begin{equation}
\sup_{\tau \in [t-s,t]} \nrm{(\theta - \theta_r)(\cdot, \tau)}_q = \sup_{\tau \in [t-s,t]} \Big(\int_{\T^d} |(\theta- \theta_r)(x, \tau)|^q \,dx\Big)^{1/q} \leq Cr^{-d/q} r^{d/q+\beta} \leq C  r^{\beta}. \label{roughlq}
\end{equation}
The above estimates, $(\ref{rough1})$, and $(\ref{roughlq})$ imply
 
\begin{equation}
II\leq C_q C A^{1/p} r^{\beta -d/q}s \exp(Csr^{\beta-1-\epsilon})\label{R}
\end{equation}

Adding the contributions of the smooth part $(\ref{S})$ and the rough part $(\ref{R})$ and choosing $\epsilon < d/q$, we have for all $\tau \in [t-s,t]$ and for some $C'_q >0,$
\begin{equation}
\sup_{f_0 \in Lip(1)}\Big|\int_{\T^d} \psi(x,\tau) f_0(x) \,dx \Big| \leq r\exp\big(C'_q A^{1/p} sr^{\beta-1-d/q}\big). \label{anylipevolve}
\end{equation}
In particular, we have 
\begin{equation}
\sup_{f_0 \in Lip(1)}\Big|\int_{\T^d} \psi(x,t-s) f_0(x) \,dx \Big| \leq r\exp\big(C'_q A^{1/p} sr^{\beta-1-d/q}\big). \label{lipevolve}
\end{equation}
It follows that for $r\leq r_0$, $(\ref{lipevolve})$ is stronger than what we need for $(\ref{evolve})$ provided 
\begin{equation}
C'_q A^{1/p} r_0^{\beta-d/q -1 +\al} \leq \delta (\beta^{-1} -1).\label{r0}
\end{equation}
\section{The Decay of the $L^p$ Norm}

In this part, we show the way it decays of the $L^p$ norm on scale $r$ is stronger than what we need for $(\ref{evolve})$. Computing the derivative of the $p^{th}$ power of the $L^p$ norm of $\psi(\cdot,\tau)$ gives
\begin{equation}
\frac{d}{d\tau} \int_{\T^d} |\psi(x,\tau)|^p \,dx = p \int_{\T^d} \Psi (x,\tau) (-\Delta)^{\al/2} \psi(x,\tau) \,dx, \label{diff}
\end{equation}
where $\Psi(x,\tau) = |\psi(x,\tau)|^{p-2} \psi(x,\tau)$ (here we used the fact that the velocity was divergence free). We also have the well-known formula 
\begin{equation}
(-\Delta)^{\al/2} \psi(x,\tau) = C_{\al} \sum_{n \in \mathbb{Z}^d} \text{p.v.} \int_{\T^d} \frac{\psi(x,\tau) -\psi(y,\tau)}{|x-y-n|^{\al +d}}\,dy. \label{laplace}
\end{equation}
See \cite{CC} for a proof of $(\ref{laplace})$. If we plug $(\ref{laplace})$ into $(\ref{diff})$ and symmetrize, we see that the derivative of the $p^{th}$ power of the $L^p$ norm is 
\begin{equation}
\frac{p}{2}C_{\al} \sum_{n \in \mathbb{Z}^d} \lim_{l \rightarrow 0} \int_{D_l} \frac{ (\Psi(x,\tau) -\Psi(y,\tau))(\psi(x,\tau) -\psi(y,\tau))}{|x-y-n|^{\al +d}}\,dy\,dx,\label{purediff}
\end{equation}
where $D_l = \{(x,y)\in \T^d \times \T^d: |x-y| \geq l \}$. Notice that the integrand in $(\ref{purediff})$ is non-negative. If $r\leq r_)$, $s\leq r^{\al}$, and $\delta \beta^{-1} < \log 2$, $(\ref{anylipevolve})$ implies that for the smooth function $\eta(\cdot, \tau) = \Psi(\cdot, \tau) \ast \chi_r$ 
\begin{equation}
\int_{\T^d} \psi(x,\tau) \eta(x,\tau) \,dx < 2r\nrm{\nabla \eta}_{\infty} = 2 r \nrm{\Psi(\cdot,\tau) \ast \nabla \chi_r}_{\infty} \leq 2r \nrm{\Psi(\cdot,\tau)}_q\nrm{\nabla \chi_r}_p. \label{smallev}
\end{equation}
The choice of $\chi_r$ implies that $2r \nrm{\nabla \chi_r}_p \leq C r^{-d} r^{d/p} = C r^{-d/q},$ for some constant $C = C(\chi)$.
Notice that $\nrm{\Psi(\cdot, \tau)}_q = \nrm{\psi(\cdot,\tau)}_p^{p-1}.$ 

We may assume $\nrm{\psi(\cdot, t-s)}_p \geq A^{1/p} \frac{r^{-d/q}}{2}$ (provided $\delta(1 + d(\beta q)^{-1}) \leq \log 2$), otherwise the evolution would already be satisfied. The maximum principle implies $\nrm{\psi(\cdot,\tau)}_p \leq \nrm{\psi(\cdot,t)}_p$ for any $\tau \in [t-s,s]$, so substituting the above inequalities into $(\ref{smallev})$ yields 
\begin{equation}
\int_{\T^d} \int_{\T^d} \psi(x,\tau) \Psi(y,\tau)\chi_r(x-y) \,dy\,dx \leq  C r^{-d/q}\nrm{\psi(\cdot,\tau)}_p^{p-1}  \leq  2 CA^{-1/p}\nrm{\psi(\cdot,\tau)}_p^p . \label{minorant}
\end{equation}
The same inequality holds if $x$ and $y$ are interchanged. Let $I(x,y,\tau)$ be the numerator of the integrand in $(\ref{purediff})$, then $I(x,y,\tau) \geq 0$ by the above comments. Notice that $|x-y|^{-\al -d} \geq c' r^{-\al} \chi_r(x-y)$, for all $x,y \in \T^d$ with some constant $c'$ depending only on $\chi$. Therefore, the kernel in $(\ref{purediff})$ dominates $c' r^{-\al} \chi_r(x-y).$ Leaving only the central cell contribution ($n=0$) in $(\ref{purediff})$ and scaling $\Lambda$ by  $C_{\al} c' r^{-\al} \frac{p}{2}$ gives   

$$\frac{d}{d\tau} \int_{\T^d} |\psi(x,\tau)|^p \,dx \geq C_{\al} \frac{p}{2} \int_{\T^d} \int_{\T^d} I(x,y,\tau)|x-y|^{-\al -d}\,dx\,dy \geq  $$
$$C_{\al} c' \frac{p}{2} r^{-\al} \int_{\T^d}\int_{\T^d}  I(x,y,\tau)  \chi_r(x-y) \,dx\,dy $$
$I(x,y,\tau)$ is a sum of four terms: two of the form $\Psi(x,\tau) \psi(x, \tau) = |\psi(x,\tau)|^p$ and two of the form $-\Psi(x,\tau) \psi(y,\tau).$ Since $\chi_r$ has mean $1$, the former terms contribute $C_{\al} c' p r^{-\al} \nrm{\psi(\cdot,\tau)}_p^{p}.$ 
The latter terms contribute no less than $-2C C_{\al} c' A^{-1/p}  p r^{-\al} \nrm{\psi(\cdot,\tau)}_p^{p}$ by $(\ref{minorant})$ . Therefore, we have the lower bound
\begin{equation}
\frac{d}{d\tau} \int_{\T^d} |\psi(x,\tau)|^p \,dx \geq C_{\al}c'p(1-2CA^{-1/p})r^{-\al} \nrm{\psi(\cdot,\tau)}_p^p.\label{normbound}
\end{equation} 
Provided we choose $p$ first, we can choose $A$ large enough so that 
\begin{equation}
1-2CA^{-1/p} > 1/2\label{A}
\end{equation} and integrate the inequality $(\ref{normbound})$ to get
\begin{equation}
\int_{\T^d} |\psi(x,t-s)|^p \,dx \leq A r^{-(p-1)d}\exp(-cpsr^{-\al}), \label{normevolve}
\end{equation}
with $c  = c(\chi,\al,p) = C_{\al} c'/2$. It follows that for $r\leq r_0$ and $s \leq r^{\al},$ $(\ref{normevolve})$ is stronger than what we need for $(\ref{evolve})$ provided 
\begin{equation}
\delta \leq \text{min}\{\beta \log 2,(1 + d(\beta q)^{-1})^{-1}\log 2, (1 + d(\beta q)^{-1})^{-1} c\} \label{delta}
\end{equation} 
This proves the Class Evolution Lemma provided $\delta$ and $r_0$ are small and $A$ is large. 
\qed

\section{The Proof of Theorem 1}

If we are given a solution $\theta$ to $(\ref{SQG})$ with initial data whose mean is $\bar{\theta}_0 \neq 0$, we define $\widetilde{\theta}(x,t) = \theta(x,t) - \bar{\theta}_0$ and an operator $\widetilde{R}(\f) = R(\f + \bar{\theta}_0)$. The modified operator $\widetilde{R}$ still satisfies the assumptions of Theorem $1$. Furthermore, $\widetilde{\theta}$ is mean zero and $\theta$ solves $(\ref{SQG})$ if and only if $\widetilde{\theta}$ solves the equation
$$ \widetilde{\theta}_t = (\widetilde{R}\widetilde{\theta}\cdot \nabla)\widetilde{\theta} - (-\Delta)^{\al/2}\widetilde{\theta}.$$ 
If we can show there is a time $T$ such that every $\widetilde{\theta} \in C^{\infty}(\T^d \times [0,T])$ can be extended to a function in $C^{\infty}(\T^d \times [0,\infty))$, then the same conclusion holds for $\theta$. It follows that we may assume that $\theta_0$ has mean zero. 

Let $\al <1$, choose $\beta >1-\al$, and then choose $p>1$ so that $\beta + \al -d/q >1$ and $q = 2^n$ for some positive integer $n$. Now we select the parameters from the Class Evolution Lemma. Choose $A$ large enough so that $(\ref{A})$ is true. Now, we choose $\delta$ small enough so that $(\ref{delta})$ is true. Finally, we choose $r_0$ sufficiently small so $(\ref{r0})$ is true.  

Since the initial data $\theta_0$ has mean zero, the maximum principle implies that the $L^{q}$ norm of $\theta$ decays exponentially. More precisely, $\nrm{\psi(\cdot,\tau)}_q \leq C(\nrm{\theta_0}_{\infty}) \exp(-\tau/q)$. Indeed, the proof of Lemma $2.4$ in \cite{CC} implies for $\theta_0$ with mean zero,
$$ \frac{d}{d\tau} \nrm{\theta(\cdot,\tau)}_q^q = -q\int_{\T^d} |\theta(x,\tau)|^{q-2} \theta(x,\tau) (-\Delta)^{\frac{\al}{2}} \theta(x,\tau)\,dx \leq -\int_{\T^d} |(-\Delta)^{\frac{\al}{4}} \theta^{\frac{q}{2}}(x,\tau)|^2\,dx.$$
Since $\widehat{\theta}(\cdot,\tau) (0) = 0$, by passing to the Fourier side we see 
$$ \frac{d}{d\tau} \nrm{\theta(\cdot,\tau)}_q^q \leq - \int_{\T^d} |\theta^{q/2}(x,\tau)|^2 \,dx = - \nrm{\theta(\cdot,\tau)}_q^q.$$
This implies $\nrm{\theta(\cdot,\tau)}_q^q \leq \nrm{\theta_0}_q^q \exp(-\tau) \leq C \nrm{\theta_0}_{\infty}^q \exp(-\tau)$. 
It follows that there is a time $T_0$ (depending on $\nrm{\theta_0}_{\infty}$) such that$|\ip{\theta(\cdot,\tau), \mathcal{U}(r)}| \leq r^{\beta}$ if $\tau \geq T_0$ and $r_0 \leq r \leq 1$. Define 
$$ T_k = T_0 + \beta r_0^{\al} \sum_{j=0}^{k-1} e^{-\delta \al j}.$$ 
We now claim that if $t\geq T_k$, then $|\ip{\theta(\cdot,t), \mathcal{U}(r)}| \leq r^{\beta}$ for $ r\geq r_0 e^{-\delta k}$. This is certainly true for $k=0$ by definition of $T_0$. Suppose that the claim is true for some $k$. Let $r \in [r_0 e^{-\delta(k+1)},r_0 e^{-\delta k})$ and $s = \beta r^{\al}$.
Suppose that $t\geq T_{k+1}$ and $\psi(\cdot,t) \in \mathcal{U}(r)$. Notice that $re^{\delta} \geq r_0 e^{-\delta k}$ and $t-s \geq T_k$. The Class Evolution Lemma and $(\ref{duality})$ imply  
$$ \Big|\int_{\T^d} \theta(x,t) \psi(x,t)\,dx \Big| = \Big|\int_{\T^d} \theta(x,t-s) \psi(x,t-s)\,dx \Big| \leq e^{-\delta\beta}(r  e^{\delta})^{\beta} \leq r^{\beta}.$$
It follows that the claim is true for all $k\geq 0$. Passing to the limit, we see that $\theta$ pairs well against any $\mathcal{U}(r)$ after the moment
$$
T = T(\al, \nrm{\theta_0}_{\infty}) = \lim_{k\rightarrow \infty} T_k = T_0 + \frac{\beta r_0^{\al}}{1-\exp(-\delta \al)}.
$$
For any time $t \geq T$, we have
\begin{equation} 
\Big|\int_{\T^d} \theta(x,t) \f(x)\, dx \Big| \leq  r^{\beta},\label{endgame}
\end{equation}
for all $\f \in \mathcal{U}(r)$ and all $0< r \leq 1$. It follows that past the moment $T(\al, \nrm{\theta_0}_{\infty})$, the solution is H\"{o}lder continuous with exponent $\beta$ with H\"{o}lder norm uniformly bounded. 

We have now shown that there is a time after which we have a uniform bound on the $C^{\beta}$ norm on the solution for $\beta$ as close to $1$ as we wish. A generalization of the argument of Constantin, Cordoba, and Wu (\cite{CCW}) shows that this is sufficient to conclude that the solution is smooth past this moment.  


\section{Concluding Remarks}

The method presented above can be used to prove that a viscosity weak solution of $(\ref{SQG})$ eventually becomes smooth and therefore eventually becomes a classical solution. The terminology and following definitions come from \cite{CC}. Given $\theta_0 \in C^{\infty}(\T^d),$ a viscosity solution to the $(\ref{SQG})$ is a weak limit (in $L^2$) of solutions to 
\begin{equation}
\left\{ \begin{array}{l}
\theta_t = (R \theta \cdot \nabla) \theta - (-\Delta)^{\alpha/2} \theta + \epsilon \Delta \theta\\
\theta(x,0)= \theta_0(x)\\
\end{array}\right., \label{WEAKSQG}
\end{equation}
as $\epsilon \rightarrow 0$. The pertubation by the Laplacian and smooth initial initial data guarantee a solution, $\theta^{\epsilon}$, to $(\ref{WEAKSQG})$ is smooth for all times. The extra dissipative term doesn't affect the above estimates which allows us to conclude $\nrm{\theta^{\epsilon}(\cdot,t)}_{C^{\delta}}\leq C$ with $\delta > 1-\al$ for all $t \geq  T(\al, \nrm{\theta_0}_{\infty})$ uniformly in $\epsilon$. The results from \cite{CCW} now give the desired regularity. In particular, the above argument gives another proof of the main result in \cite{S}.





\vspace{.1in}
\textbf{Acknowledgment.} The author is grateful to Fedor Nazarov for his guidance in the completion of this project.

\noindent \textsc{Michael Dabkowski, University of Wisconsin-Madison}\\
{\em e-mail: }\textsf{\bf dabkowsk@math.wisc.edu}


\begin{thebibliography}{03}

\bibitem{CCW} Constantin, P., Cordoba, D., and Wu, J.; 
\emph{On the critical dissipative quasi-geostrophic equation}, Indiana Univ. Math. J. \textbf{50} (2001), 97-107

\bibitem{CW} Constantin, P. and Wu, J.;
\emph{Regularity of H\"{o}lder continuous solutions of the supercritical quasi-geostrophic equation},
Ann. Inst. H. Poincare Anal. Non Lineaire, \textbf{25} (2008) No. 6, 1103-1110

\bibitem{CV} Caffarelli, L. and Vasseur, A.;
\emph{Drift diffusion equations with fractional diffusion and the quasi-geostrophic equation},
arXiv:math/0608447v1 [math.AP], 17 Aug 2006;

\bibitem{S} Silvestre, L.,
\emph{Eventual regularization for the slightly supercritical quasi-geostrophic equation}
arXiv:0812.4901v2 [math.AP], 20 Sep 2009;

\bibitem{KN} Kiselev, A. and Nazarov, F.;
\emph{A variation on a theme of Caffarelli and Vassuer}
arXiv:0908.0923v2 [math.AP], 10 Aug 2009;

\bibitem{KNV} Kiselev, A., Nazarov, F., and Volberg, A.;
\emph{Global well-posedness for the critical 2D dissipative quasi-geostrophic equation}, Iventiones Math. \textbf{167} (2007), 445-453

\bibitem{St} Stein, E., \emph{Harmonic Analysis}, Princeton University Press, 1993

\bibitem{K} Katznelson, Y., \emph{An Introduction to Harmonic Analysis}, Third Edition, Cambridge University Press, 2004

\bibitem{CS} Caffarelli, L. and Silvestre, L.; \emph{An extension problem related to the fractional Laplacian}, Comm. Partial Differential Equations \textbf{32} (2007) 1245-1260;

\bibitem{CC} Cordoba, A. and Cordoba, D.;
\emph{A maximum principle applied to quasi-geostrophic equations}, Commum. Math. Phys. \textbf{249} (2004), 511-528


\end{thebibliography}
\end{document}